\documentclass[10pt]{article}
\usepackage[english]{babel}
\usepackage[cp1251]{inputenc}
\usepackage[centertags]{amsmath}
\usepackage{latexsym}
\usepackage{amssymb}
\usepackage{amsthm}
\usepackage{amsmath}
\usepackage{amsfonts}
\usepackage{ifthen}
\usepackage{graphicx}

\theoremstyle{definition}
\newtheorem{defn}{Definition}

\theoremstyle{plain}
\newtheorem{theorem}{Theorem}
\newtheorem{lemma}{Lemma}

\renewcommand{\title}[1]{\begin{center}\textbf{#1}\end{center}}
\renewcommand{\author}[2]{\begin{center}\small\medskip\textsc{#1} \\ \medskip\textit{#2}\end{center}}
\begin{document}

\title{\LARGE MEROMORPHIC FUNCTIONS WITH SEVERAL ESSENTIAL SINGULARITIES} 
\centerline{\Large A.~A.~Kondratyuk}
\centerline{\Large Lviv National University,}
\centerline {\Large Universitetska street 1, 79000, Lviv, Ukraine}
 \Large \centerline{July 8, 2008}
{\bf Key words:}\\ Annulus, meromorphic function, Nevanlinna
characteristic, index, value distribution theory. \\{\bf
Mathematics Subject Classification (2000)} \\ 30D35 \\Preprint\\
Lviv National University,\\ Faculty of Mathematics and Mechanics
  \refstepcounter{section}
\vskip20pt \centerline{Abstract}

{ A two-parameter characteristic of functions meromorphic on
annuli is introduced and an extension of the Nevanlinna value
distribution theory for such functions is proposed.}
\section*{Introduction}\label{section1}

Meromorphic functions with $m+1$ possible essential singularities $\{c_j\},$ $c_j\in\mathbb{C},$ $j=1,2,...,m,$ $c_{m+1}=\infty,$ are considered. In other words, we consider meromorphic functions in a $m$-punctured plane. The example of such a function can be given by the composition $f\circ\mathcal{R}$ of a function $f,$ transcendental meromorphic in $\mathbb{C},$ and a rational function $\mathcal{R}$ with $m+1$ distinct poles in $\overline{\mathbb{C}}.$

We introduce a $m+1$-parameter characteristic of meromorphic functions in a $m$-punctured plane and investigate the distribution of their values. The introduced characteristic possesses the properties similar to these ones of the classical Nevanlinna characteristic.


To begin consider meromorphic functions with two possible essential singularities. Up to a linear fractional transformation we have $0$ and $\infty.$ That is, we consider meromorphic functions in the plane punctured at the origin. We will approach to both singularities by the annuli. $A_{sr}=\{z:\,s<|z|<r\}.$ In other words, $\mathbb{C}\setminus\{0\}=\bigcup\limits_{0<s<r<+\infty}A_{sr}.$

It seems that consideration of annuli $A_{\frac1\tau r},\,1\le\tau,\,1\le r,$ will be more convenient.

Note, that there are essential differences between disks and annuli in the topoligal sence which are reflected in the theory of meromorphic functions. Firstly, the fundamental (Poincr\'{e}) group of a disk is trivial, while for an annulus, we have a group isomorphic to the additive group $\mathbb{Z}.$ Secondly, the group of automorphisms of the unit disk is rich. It consists of the M\"{o}bius transformations. Therefore, the Poisson integral formula is an invariant form of the Gauss mean theorem and the Poisson-Jensen formula is an invariant form of Jensen's formula. More over, all disks are conformally equivalent. The situation with annuli is another. The group of automorphisms of an annulus is poor. If $s\neq\frac1r,$ it consists of rotations only. Annuli $A_{1r}$ and $A_{1\rho}$ is conformally equivalent iff $r=\rho.$ Therefore, the theory of meromorphic functions on annuli is more complicated that this one in disks. But many results for meromorphic functions in $\mathbb{C}$ have their counterparts for these ones in a punctured plane $\mathbb{C}^*.$ For example, the Picard theorem may be reformulated in the form. "For every non-constant meromorphic function in $\mathbb{C}^*$ there is a linear fractional transformation $\omega(\mathbb{C}^*)\subset f(\mathbb{C}^*).$" If both, $\omega$ and $f,$ are identical maps, we have the coincidence.

In order to introduce a two-parameter characteristic we use the notion of index of a meromorphic function along a circle.

\refstepcounter{section}
\section*{$\bf 1^0.$ Index of $f$ along a circle}\label{section2}

\begin{lemma}{\label{lemma1}}
Let $f$ be a function meromorphic on $\{z:\,|z|=t\},$ non identical zero. Then
\begin{eqnarray}{\label{eq1}}
\displaystyle
\nu(t,f)=\frac1\pi\int\limits_{|z|=t}{\rm Im}\left(\frac{f'(z)}{f(z)}\,dz\right)
\end{eqnarray}
is an integer.
\end{lemma}

\begin{proof} Denote $\gamma(\theta)=te^{i\theta},$ $0\le\theta\le2\pi.$ If $f$ is holomorphic on the circle $\{z:\,|z|=t\}$ without zeroes on this circle then $\Gamma(\theta)=(f\circ\gamma)(\theta),$ $0\le\theta\le2\pi,$ is a closed path, and
\begin{eqnarray*}{\label{eq2}}
&\displaystyle
\frac1{2\pi i}\int\limits_{|z|=t}\frac{f'(z)}{f(z)}\,dz
=\frac1{2\pi i}\int\limits_{0}^{2\pi}\frac{f'(\gamma(\theta)))}{f(\gamma(\theta)))}\,\gamma'(\theta)\,d\theta=\nonumber&\\
&\displaystyle
=\frac1{2\pi i}\int\limits_{0}^{2\pi}\frac{\Gamma'(\theta)}{\Gamma(\theta)}\,d\theta
=\frac1{2\pi i}\int\limits_{\Gamma}\frac{d\zeta}{\zeta}={\rm Ind}_\Gamma(0).&
\end{eqnarray*}
Taking the real parts of both sides we have
\begin{eqnarray*}{\label{eq3}}
&\displaystyle
\nu(t,f)=2\,{\rm Ind}_\Gamma(0).&
\end{eqnarray*}

As index of $\zeta=0$ with respect to $\Gamma$ is an integer we obtain the needful conclusion in this case.

If $f$ has a unique simple zero  $a=te^{i\alpha},$ then $f(z)=g(z)(z-a),$
$\displaystyle \frac{f'}{f}=\frac{g'}{g}+\frac1{z-a},$ and
\begin{eqnarray*}{\label{eq4}}
&\displaystyle
{\rm Im}\left(\frac{1}{z-a}\,dz\right)={\rm Im}\left(\frac{ite^{i\theta}d\theta}{t(e^{i\theta}-e^{i\alpha})}\right)=&\\
&\displaystyle
={\rm Re}\frac{1}{1-e^{i(\alpha-\theta)}}\,d\theta=\frac12\,d\theta,&
\end{eqnarray*}
because the transformation $\displaystyle w=\frac1{1-z}$ maps the unit circle on the straight line ${\rm Re}\,w=1/2.$ Hence, in this case
\begin{eqnarray*}{\label{eq5}}
&\displaystyle
\nu(t,f)=\nu(t,g)+1.&
\end{eqnarray*}

If $f$ has a simple pole on the considered circle, then we obtain the similar equality with "-" instead of "+". The final conclusion in the general case now follows by induction.

The value $\nu(t,f)$ is said to be \textit{index of $f$ along the circle $\{z:\,|z|=t\}.$}

If a branch of $\log f$ may by determined on the circle, then $d\log f=\frac{f'}{f}\,dz.$ Thus, ${\rm Im}\left(\frac{f'}{f}\,dz\right)=d\arg f(z)$ and $\nu(t,f)$ is the increment of $\frac1\pi\,\arg f$ along the circle.

\end{proof}

\refstepcounter{section}
\section*{$\bf 2^0.$ Version 1 of Jensen's theorem}\label{section3}

\begin{lemma}{\label{lemma2}}
Let $f$ be a function meromorphic on the closure of the annulus $A_{sr}=\{z:\,s<|z|<r\},$ and non identical zero. Then
\begin{eqnarray}{\label{eq6}}
\displaystyle
\int\limits_s^r\frac{\nu(t,f)}{t}\,dt=\frac1\pi\int\limits_{0}^{2\pi}\log|f(re^{i\theta})|\,d\theta
-\frac1\pi\int\limits_{0}^{2\pi}\log|f(se^{i\theta})|\,d\theta.&
\end{eqnarray}
\end{lemma}

\begin{proof} Assume that there are neither zeroes nor poles of $f$ on the interval $\{z=te^{i\theta},\,s\le t\le r\}.$ Then a branch of $\log f$ can by determined on this interval, and we have
\begin{eqnarray*}{\label{eq7}}
\displaystyle
\log f(re^{i\theta})-\log f(se^{i\theta})=\int\limits_s^r\frac{f'(z)}{f(z)}\,dz.&
\end{eqnarray*}
This implies
\begin{eqnarray*}{\label{eq8}}
\displaystyle
\log|f(re^{i\theta})|-\log|f(se^{i\theta})|=\int\limits_s^r{\rm Re}\left(\frac{f'(te^{i\theta})}{f(te^{i\theta})}\,e^{i\theta}\right)\,dt.&
\end{eqnarray*}

The last equality is valid for all $\theta$ from $[0,2\pi]$ except for a finite number of $\theta.$ The integration over $\theta$ yields
\begin{eqnarray*}{\label{eq9}}
&\displaystyle
\frac{1}{\pi}\int\limits_0^{2\pi}\log|f(re^{i\theta})|\,d\theta
-\frac{1}{\pi}\int\limits_0^{2\pi}\log|f(se^{i\theta})|\,d\theta=&\\
&\displaystyle
=\frac1\pi\int\limits_0^{2\pi}d\theta\int\limits_s^r\frac1t
{\rm Im}\left(\frac{f'}{f}\,ite^{i\theta}\right)dt.&
\end{eqnarray*}
Using the Fubini theorem and changing the order of integration we obtain ({\ref{eq2}}).
\end{proof}

\refstepcounter{section}
\section*{$\bf 3^0.$ Version 2 of Jensen's theorem}\label{section4}

Suppose $s\le 1\le r$ and $f$ meromorphic on the closure of the annulus $A_{sr}$ and non identical zero. Put in ({\ref{eq2}}) $r=1,$ $\tau=\frac1s$ and $t=\frac 1u.$ Then
\begin{eqnarray}{\label{eq10}}
&\displaystyle
-\int\limits_{1/s}^{1}\frac{\nu(\frac1u,f)}{u}\,du=\int\limits_1^{\tau}\frac{\nu(\frac1u,f)}{u}\,du=&\nonumber\\
&\displaystyle
=\frac1\pi\int\limits_0^{2\pi}\log|f(e^{i\theta})|\,d\theta
-\frac1\pi\int\limits_0^{2\pi}\log\left|f\left(\frac1\tau\,e^{i\theta}\right)\right|\,d\theta.&
\end{eqnarray}
Putting now in ({\ref{eq2}}) $s=1$ and subtracting ({\ref{eq10}}) we have
\begin{eqnarray}{\label{eq11}}
&\displaystyle
\int\limits_{1}^{r}\frac{\nu(t,f)}{t}\,dt
-\int\limits_1^{\tau}\frac{\nu(\frac1t,f)}{t}\,dt
=\frac1\pi\int\limits_0^{2\pi}\log|f(re^{i\theta})|\,d\theta+&\nonumber\\
&\displaystyle
+\frac1\pi\int\limits_0^{2\pi}\log\left|f\left(\frac1\tau\,e^{i\theta}\right)\right|\,d\theta
-\frac2\pi\int\limits_0^{2\pi}\log\left|f\left(e^{i\theta}\right)\right|\,d\theta.&
\end{eqnarray}

In order to obtain Version 2 of Jensen's theorem we are going to connect the notion of index $\nu(t,f)$ with the counting functions of zeroes and poles of $f.$

Let $\mathbb{T}$ be the unit circle and $n(s,r;f)$ be the number of poles of $f$ in $A_{sr}.$ It follows from the argument principle that
\begin{eqnarray}{\label{eq12}}
&\displaystyle
\nu(t,f)-\nu(1,f)=&\nonumber\\
&\displaystyle
=2\,n(1,t;\frac1f)+n(\mathbb{T},\frac1f)-2\,n(1,t;f)-n(\mathbb{T},f),
\end{eqnarray}
and
\begin{eqnarray}{\label{eq13}}
&\displaystyle
\nu(1,f)-\nu(\frac1t,f)=&\nonumber\\
&\displaystyle
=2\,n(\frac1t,1;\frac1f)+n(\mathbb{T},\frac1f)-2\,n(\frac1t,1;f)-n(\mathbb{T},f)
\end{eqnarray}
under the assumption that neither zeroes nor poles of $f$ lie on the circles $\{z:\,|z|=t\}$ and $\{z:\,z=1/t\}.$

Define
\begin{eqnarray*}{\label{eq14}}
&\displaystyle
N(\tau,r;f)=\int\limits_1^\tau\frac{n(\frac1t,1;f)}{t}\,dt+\int\limits_1^r\frac{n(1,t;f)}{t}\,dt+n(\mathbb{T},f)\log\sqrt{\tau r}.
\end{eqnarray*}

Relations ({\ref{eq11}}) -- ({\ref{eq13}}) yield the following lemma.

\begin{lemma}{\label{lemma3}}
Let $f$ be a meromorphic function on the closure of $A_{\frac1\tau r}$ and non identical zero. Then
\begin{eqnarray}{\label{eq15}}
&\displaystyle
N(\tau,r;f)-N(\tau,r;f)
=\frac1{2\pi}\int\limits_0^{2\pi}\log\left|f\left(\frac{e^{i\theta}}{\tau}\right)\right|\,d\theta+&\nonumber\\
&\displaystyle
+\frac1{2\pi}\int\limits_0^{2\pi}\log|f\left(r{e^{i\theta}}\right)|\,d\theta
-\frac1{\pi}\int\limits_0^{2\pi}\log|f\left({e^{i\theta}}\right)|\,d\theta+&\\
&\displaystyle
+\nu(1,f)\log\sqrt{\frac{\tau}{r}}.&\nonumber
\end{eqnarray}
\end{lemma}

\refstepcounter{section}
\section*{$\bf 4^0.$ Characteristic of $f$}\label{section5}

Now we are able to introduce a two-parameter characteristic of meromorphic functions in an annulus which possesses the properties like to its classical counterpart. We follow the Cartan idea. Applying ({\ref{eq15}}) to the function
$f(z)-e^{i\varphi},$ and integrating over $\varphi$ on $[0,2\pi]$ we obtain
\begin{eqnarray}{\label{eq16}}
&\displaystyle
\frac1{2\pi}\int\limits_0^{2\pi}N\left(\tau,r;\frac{1}{f-e^{i\varphi}}\right)\,d\varphi-N(\tau,r;f)=&\nonumber\\
&\displaystyle
=\frac1{2\pi}\int\limits_0^{2\pi}\log^+\left|f\left(\frac{e^{i\theta}}{\tau}\right)\right|\,d\theta
+\frac1{2\pi}\int\limits_0^{2\pi}\log^+|f\left(r{e^{i\theta}}\right)|\,d\theta-&\\
&\displaystyle
-\frac1{\pi}\int\limits_0^{2\pi}\log^+|f\left({e^{i\theta}}\right)|\,d\theta
+\frac1{4\pi}\int\limits_0^{2\pi}\nu(1,f-e^{i\varphi})\,d\varphi\log{\frac{\tau}{r}},&\nonumber
\end{eqnarray}
of course, if the last integral exists.

We will show that and evaluate the integral.

Let $f$ be a meromorphic function on $\mathbb{T}.$

Denote $E_f^+=\left\{\mathbb{T}\ni z:\,|f(z)|>1\right\},$ $E_f^0=\left\{\mathbb{T}\ni z:\,|f(z)|=1\right\}$ and $E_f^-=\left\{\mathbb{T}\ni z:\,|f(z)|<1\right\}.$

\begin{lemma}{\label{lemma4}}
Let $f$ be a meromorphic function on the unit circle $\mathbb{T},$ $f(z)\not\equiv0.$ Then
\begin{eqnarray*}{\label{eq17}}
&\displaystyle
\frac1{4\pi}\int\limits_0^{2\pi}\nu(1,f-e^{i\varphi})\,d\varphi
=\frac{1}{2\pi}\int\limits_{E_f^+}{\rm Im}\left(\frac{f'}{f}\,dz\right)
+\frac{1}{4\pi}\int\limits_{E_f^0}{\rm Im}\left(\frac{f'}{f}\,dz\right).&
\end{eqnarray*}
\end{lemma}

The proof of Lemma \ref{lemma4} needs some auxiliary results.

\begin{lemma}{\label{lemma5}}
For each function $f$ meromorphic on $\{z:\,|z|=t\},$ $f\not\equiv0,$ and each $\zeta\in\mathbb{C}$ the relation
\begin{eqnarray}{\label{eq18}}
&\displaystyle
\nu(t,f-\zeta)=\nu(t,f)-\frac1{\pi}\int\limits_{|z|=t}{\rm Im}\left(\frac{\zeta f'}{f(\zeta-f)}\,dz\right)&
\end{eqnarray}
holds.
\end{lemma}

\begin{proof} Relation (\ref{eq18}) follows immediately from (\ref{eq1}) and the idenity
\begin{eqnarray*}{\label{eq19}}
&\displaystyle
\frac{f'}{f-\zeta}=\frac{f'}{f}\left(1-\frac{\zeta}{\zeta-f}\right).&
\end{eqnarray*}
\end{proof}

\begin{lemma}{\label{lemma6}}
Let $f$ be a meromorphic function on the unit circle $\mathbb{T},$ $f(z)\not\equiv0.$ Then
\begin{eqnarray}{\label{eq20}}
&\displaystyle
\frac{1}{4\pi}\iint\limits_{\mathbb{T}\times\mathbb{T}}{\rm Re}\left(\frac{f'}{f(\zeta-f)}\,dz\,d\zeta\right)=&\nonumber\\
&\displaystyle
=-\frac{1}{2\pi}\int\limits_{E_f^-}{\rm Im}\left(\frac{f'}{f}\,dz\right)
-\frac{1}{4\pi}\int\limits_{E_f^0}{\rm Im}\left(\frac{f'}{f}\,dz\right).&
\end{eqnarray}
\end{lemma}

\begin{proof} If $z\in\mathbb{T}\setminus E_f^0,$ i.e. $|f(z)|\neq1,$ we have
\begin{eqnarray*}{\label{eq21}}
&\displaystyle
\frac1{2\pi i}\int\limits_{\mathbb{T}}\frac{d\zeta}{\zeta-f(z)}=
\begin{cases}
0,   & {\text if } \ \ |f(z)|>1,\\
1,   & {\text if } \ \ |f(z)|<1.
\end{cases}&
\end{eqnarray*}
This is true for each $z$ from $\mathbb{T}\setminus E_f^0$ except for a finite number of $z.$ Multiplying the last equality by $i \frac{f' dz}{f}$ and taking the real parts we obtain
\begin{eqnarray*}{\label{eq22}}
&\displaystyle
\frac1{2\pi i}\int\limits_{\mathbb{T}_\zeta}{\rm Re}\left(\frac{f'}{f(\zeta-f)}\,d\zeta\,dz\right)=
\begin{cases}
0,   & {\text if } \ \ z\in E_f^+,\\
-{\rm Im}\left(\frac{f'}{f}\,dz\right),   & {\text if } \ \ z\in E_f^-.
\end{cases}&
\end{eqnarray*}
for each $z$ from $\mathbb{T}\setminus E_f^0$ except for a finite number of $z.$ Therefore,
\begin{eqnarray}{\label{eq23}}
&\displaystyle
\frac{1}{4\pi^2}\int\limits_{\mathbb{T}_\zeta}\int\limits_{\mathbb{T}_z\setminus E_f^0}{\rm Re}\left(\frac{f'}{f(\zeta-f)}\,dz\,d\zeta\right)=&\\
&\displaystyle
=\frac{1}{2\pi}\int\limits_{\mathbb{T}_z\setminus E_f^0}
\frac{1}{2\pi}\int\limits_{\mathbb{T}_\zeta}{\rm Re}\left(\frac{f'}{f(\zeta-f)}\,dz\,d\zeta\right)
=-\frac{1}{2\pi}\int\limits_{E_f^-}{\rm Im}\left(\frac{f'}{f}\,dz\right).\nonumber&
\end{eqnarray}
If $|f(z)|=1,$ i.e. $z\in E_f^0,$ then $f(z)=e^{i\alpha(\theta)}.$ We will use the identity
\begin{eqnarray}{\label{eq24}}
&\displaystyle
{\rm Re}\left(\frac{f'}{f(\zeta-f)}\,dz\,d\zeta\right)
={\rm Re}\left(\frac{f'}{f}\,dz\right) {\rm Re}\left(\frac{d\zeta}{\zeta-f}\right)-&\nonumber\\
&\displaystyle
-{\rm Im}\left(\frac{f'}{f}\,dz\right) {\rm Im}\left(\frac{d\zeta}{\zeta-f}\right)&
\end{eqnarray}

If $E_f^0$ contains an open arc then $d\log|f(z)|={\rm Re}\left(\frac{f'}{f}\,dz\right)=0$ on this arc, because $\log|f(z)|=0.$ By continuity
\begin{eqnarray*}{\label{eq25}}
&\displaystyle
{\rm Re}\left(\frac{f'}{f}\,dz\right)=0&
\end{eqnarray*}
on the closure of the arc. Thus, the intersection of $E_f^0$ with the set $\left\{z:\,{\rm Re}\left(\frac{f'}{f}\,dz\right)\neq0\right\}$ is empty or consists of isolated points. Besides this,
\begin{eqnarray*}{\label{eq26}}
&\displaystyle
{\rm Re}\frac{d\zeta}{\zeta-f}={\rm Im}\frac{ie^{i\varphi}}{e^{i\varphi}-e^{i\alpha(\theta)}}
={\rm Re}\frac1{1-e^{i(\alpha(\theta)-\varphi)}}=\frac12.&
\end{eqnarray*}
Hence, relation (\ref{eq24}) implies
\begin{eqnarray}{\label{eq27}}
&\displaystyle
\frac{1}{4\pi^2}\int\limits_{\mathbb{T}\times E_f^0}{\rm Re}\left(\frac{f'}{f(\zeta-f)}\,dz\,d\zeta\right)
=-\frac{1}{4\pi}\int\limits_{E_f^0}{\rm Im}\left(\frac{f'}{f}\,dz\right).&
\end{eqnarray}
Relations (\ref{eq23}) and (\ref{eq27}) yield (\ref{eq20}).
\end{proof}

\begin{proof}[Proof of Lemma \ref{lemma4}] Applying Lemma \ref{lemma5} we have
\begin{eqnarray}{\label{eq28}}
&\displaystyle
\frac{1}{4\pi}\int\limits_{0}^{2\pi}\nu(1,f-e^{i\varphi})\,d\varphi=\frac12\nu(1,f)-&\nonumber\\
&\displaystyle
-\frac{1}{4\pi^2}\int\limits_0^{2\pi}\int\limits_{\mathbb{T}}{\rm Im}\left(\frac{f' e^{i\varphi}dz}{f(e^{i\varphi}-f)}\right)\,d\varphi.&
\end{eqnarray}
As
\begin{eqnarray*}{\label{eq29}}
&\displaystyle
d\varphi=\frac{-id\zeta}{\zeta},\qquad \zeta=e^{i\varphi},
\end{eqnarray*}
then
\begin{eqnarray*}{\label{eq30}}
&\displaystyle
-\frac{1}{4\pi^2}\int\limits_{0}^{2\pi}\int\limits_{\mathbb{T}}{\rm Im}\left(\frac{\zeta f'}{f(\zeta-f)}\,dz\,d\varphi\right)=&\nonumber\\
&\displaystyle
=\frac{1}{4\pi^2}\iint\limits_{\mathbb{T}\times\mathbb{T}}{\rm Re}\left(\frac{f'}{f(\zeta-f)}\,dz\,d\zeta\right).&
\end{eqnarray*}
Applying Lemma \ref{lemma6}, which evaluates the last integral, and (\ref{eq1}) we obtain the conclusion of Lemma \ref{lemma4} from (\ref{eq28}).
\end{proof}

Denote
\begin{eqnarray*}{\label{eq40}}
&\displaystyle
m(\tau,r;f)=m\left(\frac1\tau,f\right)+m(r,f)-2\,m(1,f),&
\end{eqnarray*}
where $m(r,f)$ is usual Nevanlinna's notation.

\begin{defn}{\label{defn1}}
Let $f$ be a meromorphic function on the closure of $A_{\frac1\tau r},$ $f\not\equiv0.$ The function
\begin{eqnarray}{\label{eq31}}
&\displaystyle
T(\tau,r;f)=N(\tau,r;f)+m(\tau,r;f)+&\nonumber\\
&\displaystyle
+c_f\log\frac{\tau}{r}\quad \tau\ge1,\quad r\ge1,.&
\end{eqnarray}
where
\begin{eqnarray}{\label{eq32}}
&\displaystyle
c_f=\frac1{2\pi}\int\limits_{E_f^+}{\rm Im}\left(\frac{f'}{f}\,dz\right)
+\frac1{4\pi}\int\limits_{E_f^0}{\rm Im}\left(\frac{f'}{f}\,dz\right)&
\end{eqnarray}
is called the \textit{characteristic of $f.$}
\end{defn}

A strange example. Let $m\in\mathbb{N}.$ Then
\begin{eqnarray*}{\label{eq33}}
&\displaystyle
T(\tau,r;z^m)=T(\tau,r;z^{-m})=m\,\frac{\log\tau+\log r}{2}.&
\end{eqnarray*}

Note, that the introduced characteristic fits not only for meromorphic functions in the punctured plane but in arbitrary annuli including punctured disks.

\begin{theorem}{\label{theorem1}}
Let $f,$ $f\not\equiv0,$ be a meromorphic function on the annulus  $A_{s_0r_0},$ $s_0<1<r_0.$ Then the function
$T(\tau,r;f)$ is non-negative, continuous, non-decreasing and convex with respect to logarithm of each variable for $1\le \tau<1/s_0,$ $1\le r<r_0,$
\begin{eqnarray*}{\label{eq34}}
&\displaystyle
T(1,1;f)=0 \ \ {\mbox and}\ \ T(\tau,r;\frac1f)=T(\tau,r;f).&
\end{eqnarray*}

If $A_{s_0r_0}=\mathbb{C}\setminus\{0\},$ the function $f$ has the meromorphic continu\-ation at the origin and $T(r,f)$ is its Nevanlinna characteristic, then
\begin{eqnarray}{\label{eq35}}
&\displaystyle
T(r,f)-2\,T(1,f)\le T(r,r;f)\le T(r,f),\quad r\ge1.&
\end{eqnarray}
\end{theorem}

\begin{proof} The indicated properties follow immediately from (\ref{eq16}), Lemma~\ref{lemma4} and Lemma~\ref{lemma3}. Moreover, relation (\ref{eq16}) implies
\begin{eqnarray*}{\label{eq36}}
&\displaystyle
\tau\,\frac{\partial T}{\partial\tau}
=\frac1{2\pi}\int\limits_0^{2\pi}n\left(\frac1\tau,1;\frac1{f-e^{i\varphi}}\right)\,d\varphi
+\frac1{4\pi}\int\limits_0^{2\pi}n\left(\mathbb{T},\frac1{f-e^{i\varphi}}\right)\,d\varphi&\\
\end{eqnarray*}
and
\begin{eqnarray*}{\label{eq37}}
&\displaystyle
r\,\frac{\partial T}{\partial r}
=\frac1{2\pi}\int\limits_0^{2\pi}n\left(1,r;\frac1{f-e^{i\varphi}}\right)\,d\varphi
+\frac1{4\pi}\int\limits_0^{2\pi}n\left(\mathbb{T},\frac1{f-e^{i\varphi}}\right)\,d\varphi&\\
\end{eqnarray*}
at the points of continuity of $n.$ Thus,
\begin{eqnarray*}{\label{eq38}}
&\displaystyle
\tau\,\frac{\partial T}{\partial\tau}+r\,\frac{\partial T}{\partial r}
=\frac1{2\pi}\int\limits_0^{2\pi}n\left(\frac1\tau,r;\frac1{f-e^{i\varphi}}\right)\,d\varphi.
\end{eqnarray*}

It seems, that the last identity with some boundary condition can be a definition of $T(\tau,r;f).$

Verify also the identity $T(\tau,r;1/f)=T(\tau,r;f).$ It follows from Version 2 of Jensen's theorem (see (\ref{eq15})) and the relation
\begin{eqnarray}{\label{eq39}}
&\displaystyle
\frac12\,\nu(1,f)=c_f-c_{\frac1f}.
\end{eqnarray}
Relation (\ref{eq35}) is in~\cite{1}.
\end{proof}

\refstepcounter{section}
\section*{$\bf 5^0.$ First Fundamental Theorem (FFT) }\label{section6}

\noindent\textbf{\textsl{FFT.}} \textit{Let $f$ be a non-constant meromorphic function in $\mathbb{C}\setminus\{0\}.$ Then for each $a\in\mathbb{C}$
\begin{eqnarray}{\label{eq41}}
&\displaystyle
N\left(\tau,r;\frac{1}{f-a}\right)+m\left(\tau,r;\frac{1}{f-a}\right)=&\nonumber\\
&\displaystyle
=T(\tau,r;f)+\varepsilon_1(\tau,r,a)+\varepsilon_2(a)\log\frac{\tau}{r},\quad\tau\ge1,\quad r\ge1&
\end{eqnarray}
where
\begin{eqnarray*}{\label{eq42}}
&\displaystyle
\left|\varepsilon_1(\tau,r,a)\right|\le4\,\log^+|a|+4\log 2,&
\end{eqnarray*}
and
\begin{eqnarray}{\label{eq43}}
&\displaystyle
\left|\varepsilon_2(a)\right|\le C,\qquad C\ {\mbox {is a constant.}} &
\end{eqnarray}}

\begin{proof} As in classical Nevanlinna theory we apply Version 2 of Jensen's theorem (see (\ref{eq15})) for $f(z)-a.$ Then
\begin{eqnarray*}{\label{eq44}}
&\displaystyle
N\left(\tau,r;\frac{1}{f-a}\right)+m\left(\tau,r;\frac{1}{f-a}\right)=N\left(\tau,r,f\right)+&\nonumber\\
&\displaystyle
+m\left(\tau,r;f-a\right)+\frac12\,\nu(1,f-a)\log\frac{\tau}{r}=&\nonumber\\
&\displaystyle
=T\left(\tau,r;f\right)+m\left(\tau,r;f-a\right)-m\left(\tau,r;f\right)+&\nonumber\\
&\displaystyle
+\left(\frac12\,\nu(1,f-a)-c_f\right)\log\frac{\tau}{r}.&
\end{eqnarray*}
Set
\begin{eqnarray*}{\label{eq45}}
&\displaystyle
\varepsilon_1(\tau,r;a)=m(\tau,r;f-a)-m(\tau,r;f)&
\end{eqnarray*}
and
\begin{eqnarray*}{\label{eq46}}
&\displaystyle
\varepsilon_2(a)=\frac12\,\nu(1,f-a)-c_f.&
\end{eqnarray*}
We have, as usually
\begin{eqnarray*}{\label{eq47}}
&\displaystyle
\left|\varepsilon_1(\tau,r,a)\right|\le4\log^+|a|+4\log2.&
\end{eqnarray*}
Applying (\ref{eq18}) for $\zeta=a$ we obtain
\begin{eqnarray*}{\label{eq48}}
&\displaystyle
\varepsilon_2(a)=\frac{1}{2\pi}\int\limits_{E_f^-}{\rm Im}\left(\frac{f'}{f}\,dz\right)
+\frac{1}{4\pi}\int\limits_{E_f^0}{\rm Im}\left(\frac{f'}{f}\,dz\right)-&\nonumber\\
&\displaystyle
-\frac{1}{2\pi}\int\limits_{\mathbb{T}}{\rm Im}\left(\frac{a f'}{f(a-f)}\,dz\right).&
\end{eqnarray*}
Since the last integral is bounded as a function of $a$ then (\ref{eq43}) is valid.
\end{proof}



\end{document}